\documentclass{amsart}[11pt]
\title{On elementary equivalence, isomorphism and isogeny of
arithmetic function fields}
\author{Pete L. Clark}
\begin{document}
\newtheorem{lemma}{Lemma}
\newtheorem{prop}[lemma]{Proposition}
\newtheorem{cor}[lemma]{Corollary}
\newtheorem{thm}[lemma]{Theorem}
\newtheorem{conj}[lemma]{Conjecture}
\newtheorem{ques}{Question}
\newtheorem{guess}[lemma]{Guess}
\newtheorem{letthm}[lemma]{Theorem}
\newtheorem{letcor}[lemma]{Corollary}
\maketitle
\newcommand{\DD}{\mathcal{D}}
\newcommand{\F}{\ensuremath{\mathbb F}}
\newcommand{\Fp}{\ensuremath{\F_p}}
\newcommand{\Fpbar}{\overline{\Fp}}
\newcommand{\PP}{\mathbb{P}}
\newcommand{\Fq}{\ensuremath{\F_q}}
\newcommand{\PPone}{\mathfrak{p}_1}
\newcommand{\PPtwo}{\mathfrak{p}_2}
\newcommand{\PPonebar}{\overline{\PPone}}
\newcommand{\N}{\ensuremath{\mathbb N}}
\newcommand{\Q}{\ensuremath{\mathbb Q}}
\newcommand{\R}{\ensuremath{\mathbb R}}
\newcommand{\Z}{\ensuremath{\mathbb Z}}
\newcommand{\SSS}{\ensuremath{\mathcal{S}}}
\newcommand{\Rn}{\ensuremath{\mathbb R^n}}
\newcommand{\Ri}{\ensuremath{\R^\infty}}
\newcommand{\C}{\ensuremath{\mathbb C}}
\newcommand{\Cn}{\ensuremath{\mathbb C^n}}
\newcommand{\Ci}{\ensuremath{\C^\infty}}
\newcommand{\U}{\ensuremath{\mathcal U}}
\newcommand{\gn}{\ensuremath{\gamma^n}}
\newcommand{\ra}{\ensuremath{\rightarrow}}
\newcommand{\fhat}{\ensuremath{\hat{f}}}
\newcommand{\ghat}{\ensuremath{\hat{g}}}
\newcommand{\hhat}{\ensuremath{\hat{h}}}
\newcommand{\covui}{\ensuremath{\{U_i\}}}
\newcommand{\covvi}{\ensuremath{\{V_i\}}}
\newcommand{\covwi}{\ensuremath{\{W_i\}}}
\newcommand{\Gt}{\ensuremath{\tilde{G}}}
\newcommand{\gt}{\ensuremath{\tilde{\gamma}}}
\newcommand{\Gtn}{\ensuremath{\tilde{G_n}}}
\newcommand{\gtn}{\ensuremath{\tilde{\gamma_n}}}
\newcommand{\gnt}{\ensuremath{\gtn}}
\newcommand{\Gnt}{\ensuremath{\Gtn}}
\newcommand{\Cpi}{\ensuremath{\C P^\infty}}
\newcommand{\Cpn}{\ensuremath{\C P^n}}
\newcommand{\lla}{\ensuremath{\longleftarrow}}
\newcommand{\lra}{\ensuremath{\longrightarrow}}
\newcommand{\Rno}{\ensuremath{\Rn_0}}
\newcommand{\dlra}{\ensuremath{\stackrel{\delta}{\lra}}}
\newcommand{\pii}{\ensuremath{\pi^{-1}}}
\newcommand{\la}{\ensuremath{\leftarrow}}
\newcommand{\gonem}{\ensuremath{\gamma_1^m}}
\newcommand{\gtwon}{\ensuremath{\gamma_2^n}}
\newcommand{\omegabar}{\ensuremath{\overline{\omega}}}
\newcommand{\dlim}{\underset{\lra}{\lim}}
\newcommand{\ilim}{\operatorname{\underset{\lla}{\lim}}}
\newcommand{\Hom}{\operatorname{Hom}}
\newcommand{\Ext}{\operatorname{Ext}}
\newcommand{\Part}{\operatorname{Part}}
\newcommand{\Ker}{\operatorname{Ker}}
\newcommand{\im}{\operatorname{im}}
\newcommand{\ord}{\operatorname{ord}}
\newcommand{\unr}{\operatorname{unr}}
\newcommand{\B}{\ensuremath{\mathcal B}}
\newcommand{\Ocr}{\ensuremath{\Omega_*}}
\newcommand{\Rcr}{\ensuremath{\Ocr \otimes \Q}}
\newcommand{\Cptwok}{\ensuremath{\C P^{2k}}}
\newcommand{\CC}{\ensuremath{\mathcal C}}
\newcommand{\gtkp}{\ensuremath{\tilde{\gamma^k_p}}}
\newcommand{\gtkn}{\ensuremath{\tilde{\gamma^k_m}}}
\newcommand{\QQ}{\ensuremath{\mathcal Q}}
\newcommand{\I}{\ensuremath{\mathcal I}}
\newcommand{\sbar}{\ensuremath{\overline{s}}}
\newcommand{\Kn}{\ensuremath{\overline{K_n}^\times}}
\newcommand{\tame}{\operatorname{tame}}
\newcommand{\Qpt}{\ensuremath{\Q_p^{\tame}}}
\newcommand{\Qpu}{\ensuremath{\Q_p^{\unr}}}
\newcommand{\scrT}{\ensuremath{\mathfrak{T}}}
\newcommand{\That}{\ensuremath{\hat{\mathfrak{T}}}}
\newcommand{\Gal}{\operatorname{Gal}}
\newcommand{\Aut}{\operatorname{Aut}}
\newcommand{\tors}{\operatorname{tors}}
\newcommand{\Zhat}{\hat{\Z}}
\newcommand{\linf}{\ensuremath{l_\infty}}
\newcommand{\Lie}{\operatorname{Lie}}
\newcommand{\GL}{\operatorname{GL}}
\newcommand{\End}{\operatorname{End}}
\newcommand{\aone}{\ensuremath{(a_1,\ldots,a_k)}}
\newcommand{\raone}{\ensuremath{r(a_1,\ldots,a_k,N)}}
\newcommand{\rtwoplus}{\ensuremath{\R^{2  +}}}
\newcommand{\rkplus}{\ensuremath{\R^{k +}}}
\newcommand{\length}{\operatorname{length}}
\newcommand{\Vol}{\operatorname{Vol}}
\newcommand{\cross}{\operatorname{cross}}
\newcommand{\GoN}{\Gamma_0( N)}
\newcommand{\GAG}{\Gamma \alpha \Gamma}
\newcommand{\GBG}{\Gamma \beta \Gamma}
\newcommand{\HGD}{H(\Gamma,\Delta)}
\newcommand{\Ga}{\Gamma^\alpha}
\newcommand{\Divo}{\Div_0}
\newcommand{\Hstar}{\cal{H}^*}
\newcommand{\txon}{\tilde{X}_0(N)}

\newcommand{\Qpi}{\ensuremath{\Q(\pi)}}
\newcommand{\Qpin}{\Q(\pi^n)}
\newcommand{\pibar}{\overline{\pi}}
\newcommand{\pbar}{\overline{p}}
\newcommand{\lcm}{\operatorname{lcm}}
\newcommand{\trace}{\operatorname{trace}}
\newcommand{\OKv}{\mathcal{O}_{K_v}}

\newcommand{\OO}{\mathcal{O}}
\newcommand{\Pic}{\operatorname{Pic}}
\newcommand{\fa}{\forall}
\newcommand{\te}{\exists}
\newcommand{\ff}{\mathfrak{f}}
\newcommand{\gk}{\mathfrak{g}_k}
\newcommand{\Gm}{\mathbb{G}_m}
\newcommand{\FPic}{\mathbf{Pic}}
\newcommand{\Spec}{\operatorname{Spec}}
\newcommand{\Res}{\operatorname{Res}}
\newcommand{\Sym}{\operatorname{Sym}}
\newcommand{\Con}{\operatorname{Con}}
\begin{abstract}
Motivated by recent work of Florian Pop, we study the connections
between three notions of equivalence of function fields:
isomorphism, elementary equivalence, and the condition that each
of a pair of fields can be embedded in the other, which we call
isogeny.  Some of our results are purely geometric: we give an isogeny 
classification of Severi-Brauer varieties and quadric surfaces.  These
results are applied to deduce new instances of ``elementary equivalence
implies isomorphism'': for all genus zero curves over a number field,
and for certain genus one curves over a number field, including some
which are not elliptic curves.
\end{abstract}
\noindent
Notation: For a field $k$, $\overline{k}$ denotes a fixed choice of separable
algebraic closure of $k$ and $\gk$ denotes the absolute Galois group of $k$.
\section{Introduction}
\noindent
\subsection{The elementary equivalence versus isomorphism problem}
A fundamental problem in arithmetic algebraic geometry is to 
classify varieties over a field $k$ up to birational equivalence, i.e., to classify
finitely generated field extensions $K/k$ up to isomorphism.  On the other 
hand,
there is the model-theoretic notion of elementary equivalence of fields -- written as
$K_1 \equiv K_2$ -- i.e., coincidence of their first-order theories.  Model theorists well 
know that elementary equivalence is considerably coarser than isomorphism:
for any infinite field $F$ there exist fields of all cardinalities elementarily
equivalent to $F$ as well as infinitely many isomorphism classes of countable 
fields elementarily equivalent to $F$.  
\\ \\
However, the fields elementarily equivalent to a given field $F$
produced by standard model-theoretic methods (Lowenheim-Skolem, ultraproducts)
tend to be rather large: e.g., any field elementarily equivalent to $\Q$
has infinite absolute transcendence degree \cite{Jensen-Lenzing}.  It is more
interesting to ask about the class of fields elementarily equivalent to a
given field and satisfying some sort of finiteness condition.  This leads us 
to the following
\begin{ques}
Let $K_1, \ K_2$ be function fields with respect to a field $k$.  Does
$K_1 \equiv K_2 \implies K_1 \cong K_2$?
\end{ques}
\noindent
On the model-theoretic side, we work in the language of fields and
\emph{not} in the language of $k$-algebras -- i.e., symbols for
the elements of $k \setminus \{0,1\}$ are not included in our alphabet.  
However, in the geometric study of function fields one certainly does
want to work in the category of $k$-algebras.  This turns out not to be 
a serious obstacle, but requires certain circumlocutions about function
fields, which are given below.  
\\ \\
By a \textbf{function field
with respect to $k$} we mean a field $K$ for which there exists
a finitely generated field homomorphism $\iota: k \ra K$ such that $k$ is 
algebraically closed in $K$, but the \emph{choice} of a particular
$\iota$ is not given.  Rather, such a choice of $\iota$ is said to give a 
\textbf{$k$-structure} on $K$, and we use the customary notation
$K/k$ to indicate a function field endowed with a particular
$k$-structure.  Suppose that $\varphi: K_1 \ra K_2$
is a field embedding of function fields with respect to $k$.  If
$k$ has the property that every field homomorphism $k \ra k$
is an isomorphism -- and fields of absolute transcendence
degree zero have this property -- then we can choose $k$-structures compatibly
on $K_1$ and $K_2$ making $\varphi$ into a morphism of $k$-algebras:
indeed, take an arbitrary $k$-structure $\iota_1: k \ra K_1$ and
\emph{define} $\iota_2 = \varphi \circ \iota_1$.   
\\ \\
Question 1 was first considered for one-dimensional function
fields over an algebraically closed base field by Duret (with subsequent related work
by Pierce), and for arbitrary function fields over a base field which is 
either algebraically closed
or a finite extension of the prime subfield (i.e., a finite field or a number
field) by Florian Pop.  They obtained the following results: 
\begin{thm}(\cite{Duret}, \cite{Pierce})
Let $k$ be an algebraically closed field, and $K_1, \ K_2$ be
one-variable function fields with respect to $k$ such that $K_1 \equiv K_2$. \\
a) If the genus of $K_1$ is different from $1$, then $K_1 \cong K_2$. \\
b) If the genus of $K_1$ is one, then so also is the genus of $K_2$,
and the corresponding elliptic curves admit two isogenies of relatively prime
degrees.
\end{thm}
\noindent
(Duret's work was formulated in the language of $k$-algebras; however,
when $k$ is algebraically closed, $k$ is definable in $K$, and the distinction
is not as critical as in the present case.)
\\ \\
The conclusion of part b) also allows us
to deduce that $K_1 \cong K_2$ in most cases, e.g. when the corresponding
elliptic curve $E_1/k$ has $\End(E_1) = \Z$.
\\ \\
The \textbf{absolute subfield} of a field $K$ is the algebraic closure
of the prime subfield ($\F_p$ or $\Q$) in $K$.  It is easy to see
that two elementarily equivalent fields must have isomorphic absolute
subfields.
\begin{thm}(\cite{Pop})
Let $K_1, \ K_2$ be two function fields
with respect to an algebraically closed field $k$ such that $K_1 \equiv K_2$.
Then: \\
a) $K_1$ and $K_2$ have the same transcendence degree over $k$. \\
b) If $K_1$ is of \textbf{general type}, $K_1 \cong K_2$.
\end{thm}
\noindent
We recall that having general type means that for a corresponding projective 
model $V/k$ with $k(V) = K$, the linear system given by a sufficiently large
positive multiple of the canonical class gives a birational embedding of
$V$ into projective space.  For curves, having general type means precisely
that the genus is at least two, so Theorem 2 does not subsume but rather
complements Theorem 1.
\\ \\
Pop obtains even stronger results (using the recent 
affirmative solution of the Milnor conjecture on $K$-theory and quadratic 
forms) in the case of finitely generated function fields.
\begin{thm}(\cite{Pop})
Let $K_1,\ K_2$ be two finitely generated fields with 
$K_1 \equiv K_2$.  Then there exist
field homomorphisms $\iota: K_1 \ra K_2$ and $\iota': K_2 \ra K_1$.  In
particular, $K_1$ and $K_2$ have the same absolute transcendence degree.
\end{thm}
\noindent
Let $K$ be a finitely generated field with absolute subfield isomorphic
to $k$.  Via a choice of $k$-structure, $K/k$ is the field of rational 
functions of a reduced, geometrically irreducible projective variety
$V/k$.  
\begin{cor}
Let $K_1/k$ be a function field of general type over either a number
field or a finite field.  Then any finitely generated field which is 
elementarily equivalent to $K_1$ is isomorphic to it.
\end{cor}
\noindent
Proof: By Theorem 3, there are field homomorphisms
$\varphi_1: K_1 \ra K_2$ and $\varphi_2: K_2 \ra K_1$, so $\Phi = 
\varphi_2 \circ \varphi_1$ gives a field homomorphism from $K_1$ to itself.  
If we 
choose $k$-structures  $\iota_i: k \hookrightarrow K_i$, on $K_1$
and $K_2$, 
then it need not be true that $\Phi$ gives a $k$-automorphism.  But since
$k$ is a finite extension of its prime subfield $k_0$, $\Aut(k/k_0)$ is 
finite, and some power of $\varphi$ induces the identity automorphism
of $k$.  In other words, there exists a dominant rational self-map
$\Phi: V/k \ra V/k$.  By a theorem of [Iitaka], when $V$ has
general type such a map must be birational.
\subsection{Isogeny of function fields}
Thus Theorem 3 is ``as good as'' Theorem 2.  But actually is is better,
in that one can immediately deduce that elementary equivalence implies
isomorphism from a weaker hypothesis than general type.    
\\ \\
Definition: We say that two fields $K_1$ and $K_2$ are \textbf{field-isogenous}
if there exist field homomorphisms $K_1 \ra K_2$ and $K_2 \ra K_1$ and denote
this relation by $K_1 \sim K_2$.  If for a field $K_1$ we have 
$K_1 \sim K_2 \implies K_1 \cong K_2$, we say $K_1$ is \textbf{field-isolated}.
If $K_1$ and $K_2$ are function fields with respect to $k$, they
are \textbf{$k$-isogenous}, denoted $K_1 \sim_k K_2$, if for some
choice of $k$-structure $\iota_1$ on $K_1$ and $\iota_2$ on $K_2$,
there exist $k$-algebra homomorphisms $\varphi_1: K_1/k \ra K_2/k$
and $\varphi_2: K_2/k \ra K_1/k$.
$K_1$ to $K_2$.  We say $K_1$ is \textbf{k-isolated} if $K_1 \sim_k K_2
\implies K_1 \cong K_2$.  Finally, if $K_1/k$ and $K_2/k$ are $k$-algebras,
we say $K_1$ is \textbf{isogenous} to $K_2$ if there exist $k$-algebra
homomorphisms $\varphi_1: K_1 \ra K_2$, $\varphi_2: K_2 \ra K_1$.
\\ \\
The distinction between field-isogeny and $k$-isogeny is a slightly 
unpleasant technicality.  It is really the notion of isogeny
of $k$-algebras which is the most natural (i.e., the most geometric),
whereas for the problem of elementary equivalence versus isomorphism,
Theorem 3 gives us field-isogeny.  There are several ways around
this dichotomy.  The most extreme is to restrict attention to 
base fields $k$ without nontrivial automorphisms, the so-called
\textbf{rigid fields}.  These include $\F_p, \ \R,\  \Q_p, \ \Q$
and ``most'' number fields.  In this case, all $k$-structures are
unique and we get the following generalization of Corollary 5.
\begin{cor}
Let $K$ be a function field with respect to its absolute subfield $k$
and assume that $k$ is rigid.  Then if $K$ is $k$-isolated, any
finitely generated field elementarily equivalent to $K$ is
isomorphic to $K$.
\end{cor}
\noindent
The assumption of a rigid base is of course a loss of generality (which
is not necessary, as will shortly become clear), but it allows us to
concentrate on the purely geometric problem of
classifying function fields $K/k$ up to isogeny.  In particular, which function 
fields are isolated?  Which have finite isogeny classes?  
\\ \\
We make some general comments on the notion of isogeny of function fields: \\ 
\\
a) The terminology is borrowed from the arithmetic theory of abelian varieties:
indeed if $K_1,\ K_2$ are function fields of principally polarized 
abelian varieties $A_1, \ A_2$, 
then they are isogenous in the above sense if and only if there is a
surjective homomorphism of group schemes with finite kernel 
$\varphi: A_1 \ra A_2$ (the point being that in this case there is
a canonical (dual) isogeny $\varphi^{\vee}: A_2 \ra A_1$).  
\\ \\
b) By a \textbf{model} $V/k$ for a function field $K/k$, we mean a 
\emph{nonsingular} projective variety with $k(V) \cong K$.  Thus the
assertion that every function field has a model relies on resolution of
singularities, which is known at present 
for transcendence degree at most two in all characteristics (Zariski,
Abhyankar) and in arbitrary dimension in characteristic zero (Hironaka).
We must emphasize that none of our results -- with the single exception
of Proposition 2d), which is itself not used in any later result -- are
conditional on resolution of singularities.  On the one hand, we are almost
entirely concerned with function fields of dimension at most two 
\emph{and} of characteristic zero; but more fundamentally, the function
fields considered here are given \emph{a priori} as function fields
of nonsingular projective varieties.
\\ \\
We can express the notion of isogeny
of two function fields $K_1/k$ and $K_2/k$ in terms of any models $V_1$
and $V_2$ by saying that there are generically finite \emph{rational}
maps $\iota: V_1 \ra V_2$ and $\iota': V_2 \ra V_1$.  
\\ \\
As usual in geometric classification problems, the easiest way to show that 
two fields $K_1$ and $K_2$ are not isogenous is not to argue directly but 
rather to find some \emph{invariant} that distinguishes between them.  It turns
out that the isogeny invariants we use are actually field-isogeny invariants.
\newcommand{\cchar}{\operatorname{char}}
\begin{prop}
Let $k$ be a field.  The following properties of a function field $K/k$
are isogeny invariants.  Moreover, when $K$ is a function field with
respect to its absolute subfield $k$, then they are also field-isogeny 
invariants. \\
a) The transcendence degree of $K/k$. \\
b) When $k$ has characteristic zero, the Kodaira dimension of a model $V/k$ for $K$. \\
c) For any model $V/k$ of $K$, the rational points
$V(k)$ are Zariski-dense. \\
d) (assuming resolution of singularities) 
For any \emph{nonsingular} model $V/k$ of $K$, there exists a 
$k$-rational point.
\end{prop}
\noindent
Proof: Part a) follows from the basic theory of transcendence bases.  As
for part b), the first thing to say is that it is \emph{false} in characteristic
$2$: there are unirational K3 surfaces [Bombieri-Mumford].  However
in characteristic zero, if $X \ra Y$ is a generically finite rational
map of algebraic varieties, then the Kodaira dimension of $Y$ is at most
the Kodaira dimension of $X$.  Moreover, the Kodaira dimension is independent
of the choice of $k$-structure.  For part c), If $X \ra Y$ is a generically finite
rational map of $k$-varieties and the rational points on $X$ are
Zariski-dense, then so too are the rational points on $Y$, so the
Zariski-density of the rational points is an isogeny invariant.
Moreover, if $\sigma$ is an automorphism of $k$, then the natural 
map $V \ra V^{\sigma} = V \times_{\sigma} k$ is an isomorphism
of abstract schemes which induces a continuous bijection
$V(k) \ra V^{\sigma}(k)$.  It follows that the Zariski-density
of the rational points is independent of the choice of $k$-structure.
For the last part, we recall the theorem of [Nishimura] that the
condition of having a simple $k$-rational point is a birational
invariant of complete varieties; moreover, as above, this condition is
independent of the choice of $k$-structure, so it suffices to consider
$K_1/k \ra K_2/k$ an embedding of function fields over $k$.
We can choose complete, \emph{normal} models $V_1/k$ and $V_2/k$
for $K_1$ and $K_2$ such that $V_2 \ra V_1$ is a morphism of varieties,
but unfortunately a simple $k$-rational point $P$ on $V_2$ could map
to a nonsimple $k$-rational point on $V_1$.  However, assuming resolution
of singularities, let $V_1/k$ be a nonsingular projective model
of $K_1/k$ and let $V_2$ be the normalization of $V_1$ in $K_2$,
so $V_1 \ra V_2$ is a morphism of $k$-varieties.  By our assumption and
by Nishimura's theorem, $V_1$ has a $k$-rational point, which maps to 
a $k$-rational point on $V_2$. 
\\ \\
The ``invariants'' of Proposition 2 are really only useful 
in analyzing the isogeny classes of varieties $V/k$ without $k$-rational
points.  For instance, two elliptic function fields $\Q(E_1)$ and
$\Q(E_2)$ have the same invariants a), b), c), d) if and only
if the groups $E_1(\Q)$ and $E_2(\Q)$ are both finite or both infinite:
this is a feeble way to try to show that two elliptic curves are not
isogenous!  
\subsection{The Brauer kernel}
In addition to the isogeny invariants of the previous subsection, we introduce
another class of invariants of a function field $k(V)$, \emph{a priori} trivial
if $V(k) \neq \emptyset$, and having the advantage that their elementary
nature is evident (rather than relying on the recent proof of the Milnor
conjecture): the Brauer kernel.
\\ \\
Let $V/k$ be a (complete nonsingular, geometrically irreducible, as always)
variety over any field $k$, and recall the exact sequence
\begin{equation}
0 \ra \Pic(V) \ra \FPic(V/k)(k) \stackrel{\alpha}{\ra} Br(k) 
\stackrel{\beta}{\ra} Br(k(V))
\end{equation}
where $\Pic(V)$ denotes the Picard group of line bundles on the $k$-scheme
$V$ and $\FPic(V/k)$ denotes the group scheme representing the sheafified
Picard group, so that in particular $\FPic(V)(k) = 
\Pic(V/\overline{k})^{\gk}$ gives the group of geometric line bundles which
are linearly equivalent to each of their Galois conjugates.  The map
$\alpha$ gives the obstruction to a $k$-rational divisor class coming from
a $k$-rational divisor, which lies in the Brauer group of $k$.  One way
to derive (1) is from the Leray spectral sequence associated
to the \'etale sheaf $\Gm$ and the morphism of \'etale sites induced by
the structure map $V \ra \Spec k$.  For details on this, see [Grothendieck].
\\ \\
\newcommand{\image}{\operatorname{image}}
We denote by $\kappa = \ker(\beta) = \image(\alpha)$ the 
\textbf{Brauer kernel} of $V$.  Some of its useful properties are:
since a $k$-rational point on $V$ defines a splitting of $\beta$,
$V(k) \neq \emptyset$ implies $\kappa = 0$.  Moreover, since it
is defined in terms of the function field $k(V)$, it is a birational
invariant of $V$.  The subgroup $\kappa$
depends on the $k$-structure on $k(V)$ as follows: if $\sigma$ is an 
automorphism
of $k$, then the Brauer kernel of $V^{\sigma} = V \times_{\sigma} k$
is $\sigma(\kappa)$.  If $k$ is a finite field, $\kappa = 0$ (since
$Br(k) = 0$).  
\\ \\
If $k$ is a number field, then $\kappa$ is a finite
group, being an image of the finitely generated group $\FPic(V)(k)$
in the torsion group $Br(k)$.  Moreover the Galois conjugacy
class of $\kappa \subset Br(k)$ is an elementary invariant of $K = k(V)$:
knowing the conjugacy class of $\kappa$ is equivalent to knowing
which finite-dimensional simple $k$-algebras $B$ (up to conjugacy) become 
isomorphic
to matrix algebras in $K$.  But if $[B:k] = n$, $B \otimes_k K$
can be interpreted in $K$ (up to $\gk$-conjugacy) via a choice
of a $k$-basis $b_1,\ldots,b_n$ of $B$ and $n^3$ structure constants
$c_{ij}^l \in k$ coming from the equations $b_i \cdot b_j = 
\sum_{l=1}^n c_{ij}^l b_l$ and the $c_{ij}^l$ themselves represented in
terms of the minimal polynomial for a generator of $k/\Q$.  Then we 
can write down the statement that $B \otimes_{k} K \cong M_n(K)$ as the
existence of an $n^2 \times n^2$ matrix $A$ with nonzero determinant and
such that 
$A(b_i \cdot b_j) = A(b_i) \cdot A(b_j)$ for all $1 \leq i,\ j \leq n$.\footnote{In
fact one can see that the conjugacy class of the Brauer kernel is an elementary
invariant whenever $k$ is merely algebraic over its prime subfield.}
\\ \\
Moreover, for any finite extension $l/k$, the conjugacy class of 
the Brauer kernel of
$V/l$ (which can be nontrivial even when $\kappa(V/k) = 0$) is again
an elementary invariant of $k(V)$.
\\ \\
If $k(V_1) \ra k(V_2)$ is an embedding of function fields, then clearly
$\kappa(V_1) \subset \kappa(V_2)$.  It follows that the Brauer kernel
is an isogeny invariant, and the Galois-conjugacy class of the Brauer
kernel is a field-isogeny invariant.
\subsection{Statement of results}
\noindent
We begin with a result relating isomorphism, isogeny, Brauer kernels
and elementary equivalence of function fields of certain geometrically
rational varieties.
\begin{thm}
For any field $k$ and any positive integer $n$, let $SB_n$ be the set of
function fields of Severi-Brauer varieties of dimension $n$ over $k$
and $Q_n$ the class of function fields of quadric hypersurfaces of
dimension $n$ over $k$. \\
a) Let $K_1, \ K_2 \in SB_n$ be cyclic elements.\footnote{A Severi-Brauer
variety $X/k$ is said to be cyclic if its corresponding division algebra $D/k$
has a maximal commutative subfield $l$ such that $l/k$ is a cyclic (Galois) extension.}
The following are 
equivalent: \\
i) $K_1 \cong K_2$. \\
ii) $K_1$ and $K_2$ are isogenous (as $k$-algebras). \\
iii) $K_1$ and $K_2$ have equal Brauer kernels. \\
b) If $K_1, \ K_2 \in Q_n$, $n \leq 2$ and the characteristic of
$k$ is not two, the following are equivalent: \\
i) $K_1 \cong K_2$. \\
ii) $K_1$ and $K_2$ are isogenous (as $k$-algebras). \\
iii) $K_1$ and $K_2$ have equal Brauer kernels, and for
every quadratic extension $l/k$, $lK_1$ and $lK_2$ have
equal Brauer kernels. \\
c) Let $K_1 \in SB_n$ and $K_2 \in Q_n, n > 1$.  Assume the characteristic
of $k$ is not two.  The following are equivalent:\\
i) $K_1 \cong K_2 \cong k(t_1,\ldots,t_n)$ are rational function fields. \\
ii) $K_1 \cong K_2$. \\
iii) $K_1$ and $K_2$ are isogenous. \\
iv) $K_1$ and $K_2$ have equal Brauer kernels. \\
\end{thm}
\noindent
\begin{cor}
Suppose $k$ is algebraic over its prime subfield.  Let $K_1 \equiv K_2$
be two function fields satisfying the hypotheses of part a), part b) or
part c) of the theorem.  Then $K_1 \cong K_2$.
\end{cor}
\noindent
Proof of the corollary: By the discussion of Section 1.3, the elementary
equivalence of $K_1$ and $K_2$ imply that their Brauer kernels are
Galois conjugate.  It follows that for any choice of $k$-structure on $K_1$, there 
exists a unique $k$-structure on $K_2$ such that we have $\kappa(K_1/k) = 
\kappa(K_2/k)$.  The theorem then implies that $K_1 \cong_k K_2$ as $k$-algebras
with this choice of $k$-structure; \emph{a fortiori} they are isomorphic as abstract
fields. \\ \\
Remarks: \\
$\bullet$ When $n = 1$, $SB_1 = Q_1$ and this class can be described equally
well in terms of genus zero curves, quaternion algebras and ternary quadratic
forms.  The essential content of the theorem when $n = 1$ is that
the Brauer kernel of a genus zero curve which is not $\mathbb{P}^1$ is cyclic of order 
two, generated by the corresponding quaternion algebra (Proposition 15).  This 
fundamental result was first proved by [Witt]. \\ \\
$\bullet$ It is well-known that the cyclicity hypothesis is satisfied for all elements of the Brauer group of a 
local or global field and for any field when $n \leq 2$.  The hypothesis
that $K_1$ corresponds to a cyclic algebra can be weakened: it
is enough to consider Severi-Brauer varieties corresponding to a solvable
crossed-product algebra [Roquette] or of various small 
degrees [Krashen].
Assuming a conjecture of Amitsur -- see Theorem 16c) -- part a) of the
theorem is valid for all Severi-Brauer function fields.
\\ \\
$\bullet$ The equivalence of bi) and bii) was first shown by [Ohm], using
results of Cassels-Pfister and Wadsworth from the algebraic theory of
quadratic forms.  (Ohm shows more, giving necessary and sufficient
conditions for one element of $Q_2$ to embed in another.)  The present author independently 
found a similar proof.  On the other hand, by
a Galois-cohomological study of the twisted forms of 
$\mathbb{P}^1 \times \mathbb{P}^1$ we give geometric proofs of
these two theorems and deduce also part iii) concerning Brauer kernels.
We emphasize that since the conjugacy class of the Brauer kernel is 
``obviously'' an elementary invariant, the instances of elementary
equivalence implies isomorphism stated in Theorem 8 are independent
of Theorem 3 (and in particular of the Milnor conjecture).
\\ \\
$\bullet$
The equivalence of i) and ii) in part b) is known also
for $n = 3$ by combining work of [Ahmad-Ohm] and [Hoffmann] -- see
[Ohm] -- but apparently not for all quadric hypersurfaces in any higher 
dimension, 
even over the rational numbers.  Condition iii) in b) is certainly false
when $n \geq 3$: the Brauer kernels of such quadrics are zero, a fact
which is used in part c).  (I owe this simple but useful observation
to Ambrus Pal.)
\begin{cor}
Let $k$ be a number field and $K$ a genus zero, one-variable
function field with respect to $k$.  Then any finitely generated field
elementarily equivalent to $K$ is isomorphic to $K$.
\end{cor}
\noindent
Proof of the Corollary: Let $L$ be a finitely generated field
such that $L \equiv K$.  Theorem $C$ applies to show that there exist
field embeddings $\iota_1: L \hookrightarrow K$ and $\iota_2: K \hookrightarrow L$.
By an appropriate choice of $k$-structures, we may view $\iota_2$ as a 
$k$-algebra morphism, hence corresponding to a morphism of algebraic curves
$C_K \ra C_L$.  By Riemann-Hurwitz, $C_L$ has genus zero, so the result follows from
Corollary 4.
\\ \\
Remark: The proof of Corollary 5 is valid for all fields $k$ finitely generated
over their prime subfield, i.e., it works also when $k$ is a finite field.
However, every genus zero curve over a finite field is isomorphic to $\PP^1$,
so in this case the result follows immmediately from Theorem 3.
\\ \\
Unfortunately the proof of Corollary 5 does not carry over
to higher-dimensional rational function fields.
when $n=1$.  Indeed, consider the case of
$K = k(t_1,\ldots,t_n) = k(\mathbb{P}^n)$, a rational function
field.  Then the isogeny class of $K$ is precisely the class of
$n$-variable function fields which are unirational over $k$.
When $n = 1$ every $k$-unirational function field is
$k$-rational, as is clear from the Riemann-Hurwitz formula and the
proof of Corollary 9 (and is well known in any case, ``Luroth's theorem'').  
If $k$ is algebraically closed of characteristic
zero, then $k$-unirational surfaces are $k$-rational, an often-noted
conseqeuence of the classification of complex algebraic surfaces
[Hartshorne, V.2.6.1].  However, for most non-algebraically closed fields
this is false, as follows from work of Segre and Manin.  Indeed,
let $K = k(S)$ be the function field of a cubic hypersurface in $\mathbb{P}^3$.
Then $K$ is unirational over $k$ if and only if for any model $S$,
$S(k) \neq \emptyset$ ([Manin 12.11]; recall that all our varieties are
smooth).  So for all $a \in k^{\times}$, the cubic surface
\[S_a: x_0^3+x_1^3+x_2^3+ax_3^3 = 0 \]
is unirational over $k$.  Segre showed that $S_a$ is $k$-rational
if and only if $a \in k^{\times 3}$; this was sharpened considerably
by [Manin, p. 184] to: $k(S_a) \cong k(S_b)$ if and only if
$a/b \in k^{\times 3}$.  Thus for any field in which 
the group of cube classes $k^{\times}/k^{\times 3}$ is infinite,
the isogeny class of $k(\mathbb{P}^2)$ is infinite.
\\ \\
As mentioned above, the isogeny invariants we have introduced here can be useful
only in classes of varieties for which $V(k) \neq \emptyset$ implies
$k(V)$ is somehow ``trivial.''  One way to make this precise is to define
an $n$-dimensional variety $V/k$ to be \textbf{prerational} if, for all field 
extensions $l/k$, $V(l) \neq \emptyset \implies l(V) \cong l(\mathbb{P}^n)$.
Is it possible that isogenous prerational varieties must be birationally
equivalent? 
\\ \\
Among one-dimensional arithmetic function fields, Question 1 is open only
for genus one curves.  By exploring the notion of an ``isogenous pair
of genus one curves'' and adapting the argument of [Pierce] in our
arithmetic context, we are able to show that elementary equivalence
implies isomorphism for certain genus one function fields.
\begin{thm}
Let $K = k(C)$ be the function field of a genus one curve
over a number field $k$, with Jacobian elliptic curve
$J(C)$.  Suppose all of the following hold:\\
$\bullet$  $J(C)$ does not have complex multiplication over $\overline{k}$. \\
$\bullet$ Either $J(C)$ is $k$-isolated \emph{or} $J(C)(k)$ is a finite group. \\
$\bullet$ The order of $C$ in $H^1(J(C),k)$ is $1, \ 2, \ 3, 4,$ or $6$. \\
Then any finitely generated field elementarily equivalent to
$K$ is isomorphic to $K$.
\end{thm}
\noindent
It goes without saying that this result is very far from a definitive treatment of the
genus one case.  Nevertheless the theorem provides evidence, convincing at least to 
the author, that the answer to Question 1
ought to be ``yes'' for all one-variable function fields over a number field. \\ \\
\noindent
Acknowledgements: The elementary equivalence versus isomorphism problem
was the topic of a lecture series and student project led by Florian
Pop at the 2003 Arizona Winter School.  In particular Corollary 8 and
Theorem 10 address questions posed by Professor Pop.  It is a pleasure
to acknowledge stimulating conversations with Pop and many other 
mathematicians at the Arizona conference, among them Abhinav Kumar,
Janak Ramakrishnan, Bjorn Poonen and Soroosh Yazdani.  The 
relevance of Galois cohomology and Brauer groups
to the classification of quadric surfaces was suggested to me by
Ambrus Pal.  This was only one of several geometric insights he shared
with me over the 2003-2004 academic year, and I am grateful for all of them.
\section{Curves of genus zero}
\noindent
The key to the case $n =1$ in Theorem 8 is the following classical (but
still not easy) result computing the Brauer kernel of a genus zero curve.
\begin{thm}([Witt])
Let $C/k$ be a genus zero curve over an arbitrary field $k$.  The
Brauer kernel of $k(C)$ is trivial if and only if $C \cong \mathbb{P}^1$.
Otherwise $\kappa(k(C)) = \{1,B_C\}$ with $B_C$ a quaternion algebra
over $k$.  Moreover the assignment $C \mapsto B_C$ gives a bijection
from the set of isomorphism classes of genus zero curves without
$k$-rational points to the set of isomorphism classes of division
quaternion algebras over $k$.
\end{thm}
\noindent
If we grant this result of Witt, the proof of Theorem 8 for function
fields of genus zero curves follows immediately: the Brauer kernel of
a genus zero curve classifies the curve up to isomorphism (and hence
its function field up to $k$-algebra isomorphism).  Moreover, the Brauer
kernel is an isogeny invariant, so genus zero curves are isogenous if and
only if they are isomorphic.  
\\ \\
We remark that Witt's theorem gives something even a bit stronger than
the $k$-isolation of the function field of a genus zero curve: it shows
that a $C_1/k$ is a genus zero curve without $k$-rational point is
not dominated by any nonisomorphic genus zero curve.  
\\ \\
We give two ``modern'' approaches to Witt's theorem: via Severi-Brauer
varieties and via quadratic forms.  We admit that part of our goal is
expository: we want to bring out the analogy between the Brauer group (of
division algebras) and the Witt ring (of quadratic forms) of a field $k$
and especially between two beautiful theorems, of Amitsur on the Brauer
group side and of Cassels-Pfister on the Witt ring side.  
\section{Severi-Brauer varieties}
\noindent
Since the automorphism groups of $M_n(k)$ and $\mathbb{P}^{n-1}(k)$
are both $PGL_{n+1}(k)$, Galois descent gives a correspondence between twisted forms
of $M_n(k)$ -- the central simple $k$-algebras -- and twisted forms of $\mathbb{P}^{n-1}$,
the Severi-Brauer varieties of dimension $n-1$.  In particular, to each Severi-Brauer
variety $V/k$ we can associate a class $[V]$ in the Brauer group of $k$, such that
two Severi-Brauer varieties of the same dimension $V_1$ and $V_2$ are \emph{isomorphic}
over $k$ if and only if $[V_1] = [V_2] \in Br(k)$. 
\\ \\
As for the birational geometry of Severi-Brauer varieties, we have the following result.
\begin{thm}([Amitsur])
Let $V_1,\ V_2$ be two Severi-Brauer varieties of equal dimension 
over a field $k$, and for $i = 1,2$ 
let $K_i = k(V_i)$ be the corresponding function field, the so-called
\textbf{generic splitting field} of $V_i$. \\
a) The subgroup $Br(K_1/k)$ of division algebras split by $K_1$ is
generated by $[V_1]$. \\
b) It follows that if $V_1$ and $V_2$ are $k$-birational, then
$[V_1]$ and $[V_2]$ generate the same cyclic subgroup of $Br(k)$. \\
c) If the division algebra representative for $V_1$ has
a maximal commutative subfield which is a cyclic Galois extension
of $k$, then the converse holds: if $[V_1]$
and $[V_2]$ generate the same subgroup of $Br(k)$, then $V_1$ and
$V_2$ are $k$-birational.
\end{thm}
\noindent
Amitsur conjectured that the last part of this theorem should remain
valid for all division algebras.  As mentioned above, there has been
some progress on this up to the present day ([Krashen]), but the general
case remains open.
\\ \\
Proof of Theorem 8 for cyclic Severi-Brauer varieties: let
$V_1/k$ and $V_2/k$ be cyclic Severi-Brauer varieties of
dimension $n$.  By Amitsur's theorem, $\kappa(V_1) = 
\kappa(V_2)$ if and only if $k(V_1) \cong_k k(V_2)$.  As in
the one-dimensional case, it follows that each of these conditions
is equivalent to $k(V_1)$ and $k(V_2)$ being isogenous (as $k$-algebras)
and in case $k$ is the absolute subfield of $k(V_1)$ and $k(V_2)$,
to $k(V_1) \equiv k(V_2)$ as fields.
\section{Quadric hypersurfaces}
\noindent
In this section the characteristic of $k$ is different from $2$.  Our second
approach to Witt's theorem is via the quadratic form(s) associated to 
a genus zero curve.
\subsection{Background on quadratic forms}
We are going to briefly review some vocabulary and results of quadratic
forms; everything we need can be found in the wonderful books  [Lam] and 
[Scharlau].  We assume familiarity with the notions of 
anisotropic, isotropic
and hyperbolic quadratic forms, as well as with the Witt ring
$W(k)$, which plays the role of the Brauer group here: it classifies
quadratic forms up to a convenient equivalence relation so that
the equivalence classes form a group, and every element of
$W(k)$ has a unique ``smallest'' representative, an anisotropic
quadratic form.  
\\ \\
The correspondence between genus zero curves over
$k$ and quaternion algebras over a field of characteristic
different from two is easy to make explicit: to a quaternion
algebra $B/k$ we associate the \textbf{ternary quadratic form} given
by the reduced norm on the trace zero subspace (of ``pure quaternions'') of
$B$.  In coordinates, the correspondence is as follows:
\[ \left( \frac{a,b}{k} \right) = 1 \cdot k \oplus i \cdot k \oplus j \cdot k
\oplus ij \cdot k \mapsto \mathcal{C}_{a,b}: aX^2 + bY^2-abZ^2 = 0. \]
By Witt cancellation, it would amount to the same to consider the
quadratic form given by the reduced norm on all of $B$; 
this quaternary quadratic form has diagonal matrix 
$\langle 1, \ a, \ b, \ -ab \rangle$.
\\ \\
On the other hand, the equivalence class of the ternary quadratic form is
\emph{not} well-determined by the isomorphism class of the curve, for the
simple reason that we could scale the defining equation of $\mathcal{C}_{a,b}$
by any $c \in k^{\times}$, which would change the ternary quadratic
form to $\langle -ca, \ -cb,cab \rangle$.  Thus at best the \textbf{similarity
class} of the quadratic form is well-determined by the isomorphism class
of $\mathcal{C}_{a,b}$, and, as we shall see shortly, this does turn out
to be well-defined.  Recall that the \textbf{discriminant} of a quadratic
form is defined as the determinant of any associated matrix, and that this quantity is
well-defined as an element of $k^{\times}/k^{\times 2}$.  It follows
that for any form $q$ of odd rank, there is a unique form similar to
$q$ with any given discriminant $d \in k^{\times}/k^{\times 2}$.  In
particular, in odd rank each similarity class contains a unique form
with discriminant $1$, which we will call ``normalized''; this leads
us to consider the specific ternary form $q_B = \langle -a, \ -b, ab \rangle$.
Moreover, to a quadratic form $q$ of any rank we can associate
its \textbf{Witt invariant} $c(q)$, which is a quaternion algebra over
$k$.  This is almost but not quite the \textbf{Hasse invariant}
\newcommand{\rank}{\operatorname{rank}} 
\[s(\langle a_1, \ldots,a_n\rangle) = \sum_{i < j} (a_i,a_j) \in Br(k) \]
but rather a small variation, given e.g. by the following \emph{ad hoc}
modifications:\footnote{Or more canonically by the theory of
Clifford algebras; see [Lam, Ch. 5].}
\[c(q) = s(q), \ \rank(q) \equiv 1, \ 2 \pmod 8, \]
\[c(q) = s(q)+(-1,-d(q)), \ \rank(q) \equiv 3, \ 4 \pmod 8, \]
\[c(q) = s(q)+(-1,-1), \ \rank(q) \equiv 5, \ 6 \pmod 8, \]
\[c(q) = s(q)+(-1,d(q)), \ \rank(q) \equiv 7, 8 \pmod 8. \]
The principal merit of $c(q)$ over $s(q)$ is that $c(q)$ is a similarity
invariant, while $s(q)$ is not.  In any case, the reader can
check that $c(q_B) = B$.  \\ \\
As a consequence of our identification of genus zero curves with quaternion
algebras, we conclude that over any field $k$, ternary quadratic forms
up to similarity are classified by their Witt invariant, and ternary
forms up to isomorphism are classified by their Witt invariant
and their discriminant (cf. [Scharlau, Theorem 13.5]).  
\\ \\
Pfister forms: For $a_1,\ldots,a_n$, we define the \textbf{$n$-fold
Pfister form} 
\[\langle \langle a_1,\ldots,a_n \rangle \rangle = \bigotimes_{i=1}^n \langle
1, \ a_i \rangle = \perp \langle a_{i_1} \cdots a_{i_k} \rangle, \]
where the orthogonal sum extends over all $2^n$ subsets of $\{1,\ldots,n\}$.
Notice that the full norm form on $B$ is $\langle 1,\ -a \rangle \otimes
\langle 1,\ -b \rangle$, a $2$-fold Pfister form.  This is good news,
since the properties of Pfister forms are far better understood than those
of arbitrary quadratic forms.  As an important instance of this, a Pfister
form is isotropic if and only if it is hyperbolic [Scharlau, Lemma 10.4].
As $n$ increases, Pfister forms become increasingly sparse among all
rank $2^n$ quadratic forms (and, obviously, among all quadratic forms), 
but observe that a quaternary quadratic form
is similar to a Pfister form if and only if it has discriminant $1$.
\\ \\
Quadric hypersurfaces: Finally, we need to link up the algebraic theory of
quadratic forms with the geometric theory of quadric hypersurfaces, our
second higher-dimensional analogue of the genus zero curves.
\\ \\
Let $q(x_1,\ldots,x_n) = a_0x_1^2+\ldots+a_nx_n^2$ be a nondegenerate
quadratic form of rank $n \geq 3$.   Let $V_q$ be the corresponding 
hypersurface in $\mathbb{P}^n$ given by $q = 0$.  $V_q$ is geometrically
irreducible and geometrically \emph{rational}.  More precisely,
$k(V_q)$ is a $k$-rational function field if and only if $q$ is isotropic:
the ``only if'' is obvious, and the converse goes as above: if we have
a single point $p \in V_q(k)$, then we can consider the family of lines
in $\mathbb{P}^{n-1}$ passing through $p$; the generic line meets $V_q$
transversely in two points, giving a birational map from
$\mathbb{P}^{n-2}$ to $V$.  However, if $n \geq 4$ then this need not be true
for every line, i.e., $V_q$ need not be isomorphic to $\mathbb{P}^{n-2}$.  
\\ \\
Every isotropic quaternary quadratic form $q$ can
be written as $H \perp g$, where $H = \langle 1, -1 \rangle$ is the hyperbolic
plane and $g$ is an arbitrary binary quadratic form; by Witt cancellation,
the equivalence classes of $g$ parameterize the isotropic quaternary
quadratic forms up to equivalence.  Since for all $c \in k^{\times}$,
$cH \cong H$, every isotropic quaternary form $q$ is similar to
$H \perp \langle 1, -d(q) \rangle$, and we conclude that isotropic quadric
surfaces are classified up to isomorphism by their discriminant.  The unique
hyperbolic representative (with discriminant $1$) is given by the equation
$x_0^2-x_1^2 + x_2^2-x_3^2 = 0$, and on this quadric we find the lines
$L_1: [a:-a:b:-b]$ and $L_2: [a:-b:-a:b]$ with intersection the single point
$[a:-a:a:-a]$: we've shown that a hyperbolic quadric surface is isomorphic
to $\mathbb{P}^1 \times \mathbb{P}^1$.
\begin{prop}
Let $q,\ q'$ be two quadratic forms over $k$.  Then $q$ is similar
to $q'$ if and only if $V_q \cong V_{q'}$.
\end{prop}
\noindent
Proof: As above, it is clear that similar forms give rise to isomorphic
quadrics.  In rank $3$ we saw that the Witt invariant, which gives
the isomorphism class of the conic, classifies the quadratic form
up to similarity.  Since a quadric surface $V$ is a twisted form of
$\mathbb{P}^1 \times \mathbb{P}^1$, the class of the canonical 
bundle in $\Pic(V)$ is represented by $K_V = -2(e_1+e_2)$, whereas the
hyperplane class of $V \subset \mathbb{P}^3$ is represented by $e_1+e_2$.
If $\varphi: V_1 \cong V_2$ is an isomorphism of quadric surfaces,
it must pull $K_{V_2}$ back to $K_{V_1}$, which, since the Picard groups
are torsionfree, implies that $e_1+e_2$ on $V_2$ pulls back to $e_1+e_2$
on $V_1$.  That is, any isomorphism of quadrics extends to an automorphism
of $\mathbb{P}^3$.  Since $\Aut(\mathbb{P}^3) = PGL_4$, this gives a similitude
on the corresponding spaces.  In rank at least $5$, the Picard group of
$V_q$ is infinite cyclic, generated by the canonical class $K_V$.  Moreover
$-K_V$ is very ample and embeds $V$ into $\mathbb{P}^{n+1}$ as a quadric
hypersurface, so again any isomorphism of quadrics extends to an automorphism
of the ambient projective space.
\\ \\
If $q$ is a rank $n$ quadratic form, we denote by $k(q)$ the function field 
$k(V_q)$ of 
the associated quadric hypersurface.  
\\ \\
If $q/k$ is a quadratic form, we say a field extension $l/k$ is
a \textbf{field of isotropy} for $q$ if $q/l$ is isotropic, or equivalently
if $l(q)$ is a rational function field.  
\\ \\
On the other hand, we say $l/k$ is a \textbf{splitting field} for $q$
if $q/l$ is hyperbolic, i.e., if $q$ lies in the ideal $W(l/k)$ of
$W(k)$ which is the kernel of the natural restriction map $W(k) \ra W(l)$.
\\ \\
The analogy with Severi-Brauer varieties and the Brauer group is
irresistible, but things are more subtle here.  Of course the function
field $k(q)$ is a field of isotropy for $q$: every variety has (generic)
rational points over its function field.  On the other hand it is not
guaranteed that $q$ becomes hyperbolic over $k(q)$.  Indeed, this
is obviously impossible unless $q$ has even rank $n = 2m$, and then
unless $d(q) = d(\mathbb{H}^m) = (-1)^m$ -- since $k$ is algebraically
closed in $k(q)$, $d(q)/(-1)^m$ does not become a square in $k(q)$ unless
it is already a square in $k$.  On the other hand, if $q$ is (similar to)
a Pfister form, then isotropy implies hyperbolicity.  So for quaternary
quadratic forms, we've shown part a) of the following result, the
analogue of Amitsur's theorem in the Witt ring.
\begin{thm}(Cassels-Pfister)[Scharlau, Theorem 4.5.4]\\
a) An anisotropic form $q$ is similar to a Pfister form if and only if
$q \in W(k(q)/k)$. \\
b) If $q$ is similar to a Pfister form and $q'$ is an anisotropic
form, then $q' \in W(k(q)/k)$ if and only if $q' \cong g \otimes q$
for some quadratic form $g$.  In particular, $W(k(q)/k)$ is
the principal ideal of $W(k)$ generated by $q$. \\
c) Let $q'$ be any quadratic form and $q$ an anisotropic quadratic
form.  If $q \in W(k(q')/k)$, then $q$ is similar to a subform
of $q'$.
\end{thm}
\noindent
(We say that $f$ is a subform of $g$ if there exists $h$ such that
$g = f \perp h$.)
\\ \\
An immediate consequence is that if $q_1$ and $q_2$
are two anisotropic Pfister forms of equal rank such that $k(q_1)$ is a 
field of isotropy for $q_2$, then $q_1$ and $q_2$ are similar.  
Applying this to the normalized norm form of a genus zero curve,
we get our second proof of Proposition 14.
\\ \\
We end this section by collecting a few more results that will be useful
for the proof of Theorem 7b).
\begin{thm}
Let $q,q'$ be quaternary quadratic forms over $k$ with common 
discriminant $d$, and put $l = k(\sqrt{d})$. \\
a)\ [Scharlau, 2.14.2] $q$ is anisotropic if and only if $q_l$ is
anistropic. \\
b)\ [Wadsworth] If $q'/l$ is similar to $q/l$, then $q$ is similar to $q'$. \\
c)\ [Wadsworth] If $q$ is anistropic and $k(q) \cong k(q')$, then
$q$ is similar to $q'$.
\end{thm}
\subsection{An algebraic proof of Ohm's theorem}
We begin the proof of Theorem 7b) by explaining how the results
we have recalled on quadratic forms can be used to deduce the
theorem of [Ohm] on the isogeny classification of quadric surfaces.
Indeed, thanks to the remarkable Theorem 15c), the classification
result is more precise than we have let on.
\begin{thm}(Ohm)
Let $q, \ q'$ be two nondegenerate quaternary quadratic forms
over $k$ with isogenous function fields.  Then either: \\
a) $q$ and $q'$ are both isotropic, so $k(q) \cong k(q') \cong k(t_1,t_2)$, 
or\\
b) $q$ and $q'$ are both anisotropic in which case $V_q \cong V_{q'}$,
i.e., $q$ and $q'$ are \emph{similar}.
\end{thm}
\noindent
That is, except in the case when both function fields are rational,
quadric surfaces with isogenous function fields are not only
birational but \emph{isomorphic}.
\\ \\
Proof: Since isotropic quadric function fields are rational and the
condition of being isotropic (i.e., of having a $k$-rational point) is
an isogeny invariant, we need only consider the case when both
$q$ and $q'$ are anisotropic quaternary quadratic forms.  The proof
divides into further cases according to the values of the discriminants
$d = d(q), \ d' = d(q')$.
\\ \\
The first case is $d = d' = 1$ (as elements of $k^{\times}/k^{\times 2}$).
In this case $q$ and $q'$ are both similar to Pfister forms.  If they are
isogenous over $k$, \emph{a fortiori} they are isogenous over $k(q')$,
and since $q'$ becomes isotropic over $k(q')$, so does $q$.  Since $q$
is a Pfister form, this implies $q \in W(k(q')/k)$, and by Theorem 14c)
we conclude that $q$ and $q'$ are similar.
\\ \\
Suppose $d = d' \neq 1$.  Let $l = k(\sqrt{d})$.  By Theorem 15a),
$q/l$ and $q'/l$ remain anisotropic.  Moreover they are now Pfister
forms, so the previous case applies to show that $q/l$ and $q'/l$
are similar.  But now Theorem 15b) tells us that $q$ and $q'$
are already similar over $k$!
\\ \\
The last case is $d \neq d'$.  Since the discriminant is a similarity
invariant among quaternary quadratic forms, we must show that
this case cannot occur, i.e., that two anisotropic quadratic
forms with distinct discriminants cannot be isogenous.  Let
$l = k(\sqrt{d})$; it suffices to show that
$q/l$ and $q'/l$ are nonisogenous.  Again, Theorem 15a) implies
that $q/l$ remains anisotropic, whereas we may assume that
$q'/l$ is anistropic, for otherwise they could not be isogenous.
We finish as in the first case: by construction $q/l$ is 
(similar to) an anistropic Pfister form, so $q/l \in W(l(q')/l)$
and the Cassels-Pfister theorem implies that $q/l$ and $q'/l$
are similar, but their discriminants are different, a contradiction.
\section{Geometry and Galois cohomology of quadric surfaces}
\noindent
Our strategy for proving Theorem 8b) in full is in fact
to make the proof of the previous subsection geometric: that is,
we will use Brauer kernels to give proofs of Theorems 14 and 15
in the case of quaternary quadratic forms.  The fact that two-dimensional
quadric function fields are classified by their Brauer kernels over
$k$ and over all quadratic extensions of $k$ will come as a byproduct
of these proofs.
\subsection{Preliminaries on twisted forms}
The first step is to consider not just
the quadric surfaces over $k$, but the larger set of all
twisted forms of the hyperbolic surface $\mathbb{P}^1 \times \mathbb{P}^1$.
\\ \\
So let $\mathcal{T} = \mathcal{T}(\PP^1 \times \PP^1)$ be the set of all
Galois twisted forms of $\PP^1 \times \PP^1$, i.e., the set of all
varieties $X/k$ such that
$X/\overline{k} \cong \PP^1 \times \PP^1$.  
We saw in the previous section that every quadric surface $V_q$ is an element
of $\mathcal{T}$.  (More precisely, every quadric surface becomes isomorphic
to $\mathbb{P}^1 \times \mathbb{P}^1$ after an extension with
Galois group $1$, $\Z/2\Z$ or $\Z/2\Z \times \Z/2\Z$, and no anistropic
quadric surface with nontrivial discriminant splits over a quadratic
extension.)
\\ \\
By Galois descent, $\mathcal{T} = H^1(k,G)$, where $G$ is the automorphism
group of $\PP^1 \times \PP^1$.  $G$ is a semidirect product: 
\[1 \ra PGL_2^2 \ra G \ra (\Z/2\Z) \ra 1, \]
where the $PGL_2^2$ gives automorphisms of each factor separately, and a splitting
of the sequence is given by the involution of the two $\PP^1$ factors.  Thus
we have a split exact sequence of pointed sets
\[1 \ra QA(k)^2 \ra \mathcal{T} \stackrel{d}{\ra} k^{\times}/k^{\times 2}, \]
where $QA(k)$ stands for the set of all quaternion algebras over $k$, and the map
$d$ can be viewed as a discriminant map.
The splitting just means that we have an injection
$k^{\times}/k^{\times 2} \hookrightarrow \mathcal{T}$: we choose the embedding
corresponding to the subset of all isotropic quadric surfaces (we have seen that these
are parameterized by their discriminant).  
\\ \\
The part of $\mathcal{T}$ in the kernel of $d$ is easy to understand: we just
take two different twisted forms $C_1, \ C_2$ of $\PP^1$ -- i.e., two genus
zero curves over $k$ -- and put $X = C_1 \times C_2$.  Using
Witt's theorem, we can identify the Brauer kernel of such a surface: 
$\kappa(k(C_1 \times C_2)) = \langle B_{C_1}, \ B_{C_2} \rangle$.
\\ \\
For any twisted form $X$, let $N = \FPic(X)(\overline{k})$ be the Picard
group of $X/\overline{k}$ viewed as a $G_k$-module.  As abelian group,
$N$ is isomorphic to $\Pic(\PP^1 \times \PP^1) = e_1^{\Z} \oplus e_2^{\Z}$,
where $e_1$ and $e_2$ represent the two rulings.  Write $N(k) := N^{G_k}$
for the $G_k$-equivariant line bundles on $X/\overline{k}$, so $N(k)$
is a free abelian group of rank at most $2$.  The rank is at least
one, since $e_1 + e_2 \in N(k)$: the only two elements of the N\'eron-Severi
lattice with self-intersection $2$ are $\pm(e_1+e_2)$, and $e_1+e_2$ is distinguished
from $-(e_1+e_2)$ by being ample; both of these properties are preserved by the
$G_k$-action.  Moreover, since $N$ is torsion free, for
any $L \in N(\overline{k})$ and any $n \in \Z^+$, $L \in N(k) \iff nL \in N(k)$.
In particular, the rank of $N(k)$ is $2$ if and only if $N(\overline{k})$
is a trivial $G_k$-module.
\\ \\
Claim: $N(k)$ has rank $2$ if and only if $d(X) = 1$. \\ \\
Proof: If $d(X) = 1$, $X = C_1 \times C_2$, and choosing any point
$p_2 \in C_2(\overline{k})$, for any $\sigma \in G_k$, $\sigma(C_1 \times p_2) =
C_1 \times \sigma(p_2)$, so that the Galois action preserves the horizontal
ruling; the same goes for the vertical ruling.  The converse is similar: to say
that $\sigma \in G_k$ acts trivially on the class of $[e_1]$ and $[e_2]$ is to
say that it does not interchange the rulings, hence lies in the subgroup
$PGL_2^2$ of $G$.
\\ \\
Look now at the rank one case, where $N(k) = (e_1+e_2)^{\Z}$.  From
the basic exact sequence (1), it follows that the Brauer kernel of $X$
is precisely the obstruction to $e_1+e_2$ coming from
a line bundle.  Since $-2(e_1+e_2)$ is represented by the canonical bundle,
we get that for all $X \in \mathcal{T}$,  $\kappa(X) \subset Br(k)[2]$.
\\ \\
Claim: $\alpha(e_1+e_2) = 0$ if and only if $X$ is a quadric surface.
\\ \\
Proof: On $\PP^1 \times \PP^1$, $H = e_1+e_2$ is very ample and gives the
embedding into $\PP^3$ as a degree $2$ hypersurface.  It follows that
as soon as the class of $[e_1+e_2]$ is represented by a $k$-rational
divisor, the same holds $k$-rationally, i.e., $X$ is embedded in 
$\PP^3$ as a degree $2$ hypersurface.  For the converse, just
cut the quadric by a hyperplane to get a rational divisor in the
class of $e_1+e_2$.
\\ \\
Let us sum up our work on Brauer kernels of quadric surfaces.
\begin{prop}
Let $X/k$ be a quadric surface.  If $d(X) \neq 1$, then the Brauer
kernel is trivial.  If $d(X) = 1$, then $X \cong C \times C$ and
is classified up to isomorphism by its Brauer kernel $\kappa(X) =
\{1,B_C\}$.
\end{prop}
\noindent
Proof: We just need to remark that when $d(X) = 1$, $X = C_1 \times 
C_2$, and since $\alpha(e_1) = B_{C_1}, \ \alpha(e_2) = B_{C_2}$ are
$2$-torsion elements of $Br(k)$, the fact that $\alpha(e_1+e_2) = 0$
implies $\alpha(e_1) = \alpha(e_2)$, so that $C_1 \cong C_2$.
\subsection{The proof of Theorem 8b}
First we give a geometric proof of Theorem 14 for
quaternary quadratic forms: Since ``Pfister quadrics'' are just 
those isomorphic to $C \times C$,
where $C$ is a genus zero curve, we can turn our previous argument on its head
and deduce part b) of the Cassels-Pfister theorem in rank $4$
from Witt's theorem.  (We saw earlier that part a) is easy in rank 4.)
Now let $q'$ be a rank $4$ quadratic form and $q$ an anisotropic
quadratic form such that $q \in W(k(q')/q')$.  But since $k$ is algebraically
closed in $k(q')$, this implies that $d(q) = 1$, so that $q$ is (similar to)
a Pfister form, with corresponding quadric $C \times C$.  If $d(q') = 1$
also, this reduces again to Proposition 14, so assume that $d(q') = d \neq 1$
and let $l = k(\sqrt{d})$.  Consider the basic exact sequence
\[0 \ra \Pic(V_2) \ra \FPic(V_2)(k) \stackrel{\alpha}{\ra} Br(k) 
\stackrel{\beta}{\ra} Br(k(V_2)). \]
The hypothesis that $q$ splits in $k(V_2)$ means that $B_C$ is an
element of the Brauer kernel of $k(V_2)$.  But being a quadric
surface with nontrivial discriminant, $\kappa(k(V_2)) = 0$, a contradiction.
\\ \\
Proof of Theorem 15a) for quaternary forms: let $q_1,\ q_2$ 
be quaternary forms with common discriminant $d$ and corresponding quadrics
$V_1, \ V_2$; put $l = k(\sqrt{d})$.
\\ \\
First we must show that if $V_1/l$ is isotropic, then $V_1/k$ was isotropic.
But if $X(l)$ is nonempty, then since the discriminant is $1$ over $l$, then 
$X$ splits over $l$.  So 
we can choose rational curves $C_1$, $C_2$ over $l$ such that $\sigma(C_1) 
= C_2$.  
But then $\sigma(C_1 \cap C_2) = \sigma(C_1) \cap \sigma(C_2) = 
C_2 \cap C_1 = C_1 \cap C_2$ gives a $k$-rational point.   
\\ \\
We now give geometric proofs of Wadsworth's results, Theorem 15 parts b) and 
c).  The isotropic
case of Theorem 15b) is easy, since isotropic quadric surfaces are classified
by their discriminant.  Since we know that two anisotropic Pfister
quadrics are birational if and only if they are isomorphic, Theorem 15c)
follows from Theorem 15b), and we are reduced to showing the following.

\begin{prop} 
Let $V/k, W/k$ be two anisotropic quadric surfaces with common 
discriminant $d$; put $l = k(\sqrt{d})$.  If $V_1/l \cong V_2/l$, then $V_1 
\cong V_2$.
\end{prop}
\noindent
Proof: We write $\sigma$ for the nontrivial element of $G_{l/k}$.  Let
$\mathcal{S}$ be the set of all $l/k$ twisted forms of $V$, and let
$\mathcal{S}_d \subset \mathcal{S}$ be the subset of twisted forms $W$
with $d(W) = d(V)$.  We claim that $\mathcal{S}_d = {V}$, which gives the 
result we want.  (In fact it is a stronger result, since we are
\emph{a priori} allowing twisted forms which are not quadric surfaces.)
\\ \\
To prove the claim we clearly may ``replace'' $V$ by any element
of $\mathcal{S}_d$.  A convenient choice is the variety $V_1/k$ constructed
as follows: let $B = c(V)$ be the Witt invariant of the quadric $V$, and let
$C/k$ denote the corresponding genus zero curve.  Let $V_1 := \Res_{l/k}(C / l)$
be the $k$-variety obtained by viewing $C$ as a curve over $l$ and then taking the
Weil restriction from $l$ down to $k$.\footnote{It is a byproduct of the proof
that $V$ is a quadric surface.  On the other hand, if we started with a genus zero curve
$C$ whose corresponding quaternion algebra was in $Br(l) \setminus Br(l)^{G_{l/k}}$,
then the restriction of scalars construction would yield a twisted form $V_1/k$
such that $V_1/l \cong C \times C^{\sigma}$ is not a quadric surface (even)
over $l$.}
\\ \\
We have that $V_1/l \cong C \times C$.  Let $G = \Aut(V_1)$.  It is convenient
(and correct!) to view $G$ as an algebraic $k$-group scheme.  In particular
this gives the $l$-valued points $G(l)$ the structure of a $G_{l/k}$-module,
and this Galois module structure is highly relevant, since
$\mathcal{S} = H^1(l/k,G(l))$.  Indeed we have a short exact sequence of
$k$-group schemes
\begin{equation}
1 \ra K \ra G \ra \Z/2\Z = \Sym \{e_1, \ e_2 \} \ra 1
\end{equation}
obtained by letting automorphisms of $V_1$ act on the $G_{l/k}$-set 
of rulings $\{e_1, \ e_2\}$; this exact sequence is of course a twisted
analogue of the exact sequence considered in 5.1.  In particular, we
still have that $K$ is the connected component of $G$, a linear
algebraic group scheme; $K(l)$ is, as a group, isomorphic
to $\Aut(C)(l)^2 = PGL(B)(l)^2$, where the group $PGL(B)$ is the twisted
analogue of $PGL_2$ defined by the short exact sequence
\[1 \ra \Gm \ra B^{\times} \ra PGL(B) \ra 1. \]
However, the $G_{l/k}$-module structure on $K(l)$ is twisted: since $\sigma$
interchanges $e_1$ and $e_2$, it also interchanges the two factors of $PGL(B)$.
That is, $\sigma(x,y) = (\sigma(y),\sigma(x))$, so \[K(k) = \{(x,\sigma(x)) \ |
\ x \in PGL(B)(l)\} \cong PGL(B)(l), \]
and one finds that $K/k = \Res_{l/k} \Aut(C/l)$.  But then Shapiro's Lemma implies
\[\#H^1(l/k,K(l)) = \#H^1(l/l,\Aut(C/l)) = 1. \]
Taking $l$-valued points and then $G_{l/k}$-invariants in (2), one gets an exact
cohomology sequence, of which a piece is
\[H^1(l/k,K(l)) \ra G \ra H^1(l/k,\Z/2\Z) = \pm 1 \]
where the latter map is the twisted analogue of the discriminant map.  The exactness means
that $H^1(l/k,K(l)) = S_d$, so we are done.
\\ \\
End of the proof of Theorem 8b): At this point, we have given a second
proof of Ohm's Theorem 16.  It remains to see that function fields
of quadric surfaces $k(X)$ are classified by their Brauer kernels over $k$ and 
over all quadratic extensions of $k$.  Suppose $k(q)$ and $k(q')$
are \emph{non}-isomorphic function fields of quadric surfaces.  If
one is isotropic and the other is anisotropic, then the isotropic one
has trivial Brauer kernels over all extension fields of
$k$, whereas the anisotropic one has a Brauer kernel of order two over
$k(\sqrt{d})$.  So suppose that both are anisotropic.  If $d(q) \neq d(q')$,
then over $k(\sqrt{d})$, $q$ has nontrivial Brauer kernel and $q'$ has
trivial Brauer kernel.  If their discriminants are the same, then by
Proposition 18,
$l(q)$ and $l(q')$ remain nonisomorphic, so have distinct nontrivial
Brauer kernels.  This shows the equivalence of the first three conditions
in part b) of Theorem 8.  
\\ \\
Remark: With a bit more care, one should be able to show
``isogeny implies birationality'' for all $\overline{k}/k$
twisted forms of $\PP^1 \times \PP^1$.  One uses similar methods
to the above, the only new wrinkle being that there are some
pairs $V_1, V_2 \in \mathcal{T}(\PP^1 \times \PP^1)$ with the
same isogeny invariants (i.e., equal Brauer kernels over
all extensions of $k$) and which are non-isomorphic
but still birational.  For instance, in the case
$d(V_1) = d(V_2) = 1$, we have $V_1 = C_1 \times C_2$, 
$V_2 = C_3 \times C_4$ (all genus zero curves), and
setting the Brauer kernels means precisely
$\{C_1,C_2\} = \{C_3,C_4\}$.  The only problematic
case is $C \times C$ versus $C \times \PP^1$, which are
evidently non-isomorphic.  However, they are birational:
let $\pi_1: C \times C \ra C$ be projection onto the first
factor; the generic fibre of $\pi_1$ is a genus
zero curve over $k(C)$.  Since there is an obvious
section -- namely $c \mapsto (c,c)$, this curve is isomorphic
to $\PP^1$ over $k(C)$, i.e., $k(C \times C) = k(C)(t) =
k(C \times \PP^1)$.  (This elegant argument is due
to Colliot-Th\'el\`ene.)

\section{Comparing quadrics and Severi-Brauer varieties}

\subsection{The proof of Theorem 8c)}
For the proof of Theorem 8c), it suffices to show that for any
$n > 1$, if $K_1 = k(V_1)$ is the function field of a nontrivial
Severi-Brauer variety and $K_2 = k(V_2)$ is the function field
of an anistropic quadric hypersurface, then $\kappa(k(V_1)) \neq
\kappa(k(V_2))$. 
\\ \\
But recall that the Picard group of a quadric hypersurface $V_2/k$ in dimension
at least $3$ is \emph{generated} by the canonical bundle, so the natural
map $\Pic(V_2) \ra \FPic(V_2)(k)$ is an isomorphism and $\kappa(k(V_2)) = 0$.
On the other hand, a nontrivial Severi-Brauer variety has a nontrivial
Brauer kernel, the cyclic subgroup generated by the 
corresponding Brauer group element.
\\ \\
When $n = 2$, the Brauer kernel of a nontrivial Severi-Brauer surface
is cyclic of order $3$, whereas the Brauer kernel of any quadric is
$2$-torsion.  
\subsection{Brauer kernels and the index}
Earlier we mentioned the fact that if $V$ has a $k$-rational point,
$\kappa(k(V)) = 0$.  This statement can be refined in terms of 
the \textbf{index} of a variety $V/k$, which is the least positive degree
of a $\gk$-invariant zero-cycle on $V$; equivalently, it is the
greatest common divisor over all degrees of finite field extensions 
$l/k$ for which $V(l) \neq \emptyset$.  Note then that the index is
a (field-)isogeny invariant.  Suppose $l/k$ is a finite
field extension of degree $n$ such that $V(l) \neq \emptyset$.
Then $\kappa(k(V)) = Br(k(V)/k) \subset Br(l/k)$.  It follows
that the index of $V/k$ is an upper bound for the index
of any element of the Brauer kernel of $k(V)$ (recall that
the index of a Brauer group element is the square root of
the $k$-vector space dimension of the corresponding division
algebra $D/k$).  In particular varieties with a $k$-rational
zero-cycle of degree one have trivial Brauer kernel.
\\ \\
Notice that quadrics and Severi-Brauer varieties have a very special
property among all varieties: namely the existence of a rational
zero-cycle of degree one implies the existence of a rational point.
For Severi-Brauer varieties, it is part of the basic theory of
division algebras that the index of a division algebra is equal to
the greatest common divisor over all degrees of splitting fields
(and moreover the gcd is \emph{attained}, by any maximal subfield
of $D/k$).  For quadrics -- whose index is clearly at most $2$
-- this follows from Springer's theorem, that an anistropic quadratic 
form remains anisotropic over any finite field extension of odd degree.
\\ \\
To see how ``special'' this property is, observe that every variety
over a finite field has index one, since the Weil bounds (it is
enough to consider curves) imply that if $V/\F_q$ is a smooth projective
variety, $V(\F_{q^n}) \neq \emptyset$ for all $n \gg 0$, and in particular
there exists $n$ such that $V/\F_q$ has rational zero cycles of coprime
degrees $n$ and $n+1$.  This gives amusingly convoluted proofs of
the familiar facts that the Brauer group of a finite field is trivial and that every
quadratic form in at least three variables over a finite field is isotropic.
\section{Curves of Genus One}
\noindent
In this section we suppose that all fields have characteristic zero.
\subsection{Preliminaries on genus one curves}
Let $K = k(C)$ be the function field of a genus one curve.
Recall that $C$ can be given the structure of an elliptic curve
if and only if $C(k) \neq \emptyset$.  Moreover, if $C$ is an arbitrary
genus one curve, we can associate to it an elliptic curve, its
\textbf{Jacobian} $J(C) = \FPic^0(C)$, the group scheme representing
the subfunctor of $\FPic(C)$ consisting of divisor classes of degree zero.
The Riemann-Roch theorem gives a canonical identification $C = \FPic^1(C)$;
with this identification, $C$ becomes a principal homogeneous space (or torsor)
over $J(C)$.  By Galois descent, the genus one curves $C/k$ with Jacobian
a given elliptic curve $E$ are parameterized by the Galois cohomology
group $H^1(k,E)$.  There is a subtlety here: $H^1(k,E)$ parameterizes
isomorphism classes of genus one curves endowed with the structure of a 
principal homogeneous space for $E$, so a genus one curve up to isomorphism
corresponds to an orbit of $\Aut(E/k)$ on $H^1(k,E)$.  We will assume
that $\Aut(E/k) = \pm 1$ (this excludes only the notorious $j$-invariants
$0$ and $1728$) -- later we will exclude all elliptic curves with
complex multiplication over the algebraic closure of $k$.  So $[C]$ and
$-[C]$ are in general distinct classes in $H^1(K,E)$ but represent
isomorphic genus one curves.
\\ \\
If a genus one curve has a $k$-rational zero-cycle of degree one,
then by Riemann-Roch it is an elliptic curve, i.e., index one implies
the existence of rational points for genus one curves.  Another important
numerical invariant of a genus one curve $C/k$ is its \textbf{period},
which is simply the order of $[C]$ in the torsion group $H^1(k,J(C))$.
It is not hard to see that, like the index, the period is an isogeny
invariant.
\\ \\
Recall that an isogeny of elliptic curves (in the usual sense) is
just a finite morphism of varieties $\varphi: (E_1,O_1) \ra (E_2,O_2)$ preserving the distinguished 
points.  But notice that if $f: E_1 \ra E_2$ is any finite morphism 
of genus one curves with rational points, it can be viewed as an isogeny
by taking $O_2 = f(O_1)$.  Moreover, if $f: E_1 \ra E_2$ is a finite
morphism of elliptic curves, then there is an induced map
$\Pic^0(f) = \Pic^0(E_2) \ra \Pic^0(E_1)$.  Since any elliptic curve
is isomorphic to its Picard variety, this explains why our notion of
an isogenous pair of elliptic function fields is consistent with
the usual notion of isogenous elliptic curves: the morphism in the other
direction is guaranteed.
\\ \\
But if $\varphi: C_1 \ra C_2$ is a morphism of genus one curves
without rational points, then since $C_2$ is not isomorphic to
$\Pic^0(C_2)$, the existence of a finite map $\varphi': C_2 \ra C_1$
is not guaranteed.  Indeed, it need not exist: let $C$ be a genus
one curve of period $n > 1$.  Then the natural map $[n]: C =
\FPic^1(C) \ra \FPic^n(C) \cong J(C)$ gives a morphism of
degree $n^2$ from $C$ to its Jacobian.  Since $J(C)(k) \neq \emptyset$,
there is no map in the other direction.  So the classification of
genus one curves up to isogeny is more delicate than the analogous
classification of elliptic curves.  We content ourselves here with the
following result.
\begin{prop}
Let $C,\ C'/k$ be two genus one curves with common Jacobian $E$, and assume
that $E$ does not have complex multiplication over $\overline{k}$.  Then
there exists a degree $n^2$ \'etale cover $C \ra C'$ if and only if
$[C'] = \pm n[C]$ as elements of the Weil-Chatelet group $H^1(k,E)$.
\end{prop}
\noindent
Proof: As we saw above, there is a natural map $\psi_n: C = \Pic^1(C) \ra 
\Pic^n(C)$ induced by the map $D \mapsto nD$ on divisors.  Upon
basechange to the algebraic closure and up to an isomorphism, this map
can be identified with $[n]$ on $J(C)$, so it is an \'etale cover of
degree $n^2$.  Keeping in mind that $n$ could be negative, corresponding
to a twist of principal homogeneous structure by $[-1]$, we get the
first half of the result.
\\ \indent
For the converse, let $\pi: C \ra C'$ be any finite \'etale
cover.  Choosing $P \in C(\overline{k})$ and its image $P' = \pi(P) \in
C'(\overline{k})$, $\pi/\overline{k}: C/\overline{k} \ra C/\overline{k}$
is an elliptic curve endomorphism.  By assumption on $E$,
$\pi/\overline{k} = [n]$ for some integer $n$, and its kernel $E[n]$
is the unique subgroup isomorphic to $\Z/n\Z \oplus \Z/n\Z$.  It follows
that the map $\pi: C \ra C'$ factors as $C \ra C/E[n] \ra C'$, hence
$C/E[n] \ra C'$ is an isomorphism of varieties.  It is not necessarily
a morphism of principal homogeneous spaces: it will be precisely when $n > 0$
in $\pi/\overline{k}$ above.  Taking into account again the possibility
of $n < 0$ gives the stated result.
\begin{cor}
Let $C/k$ be a genus one curve with non-CM Jacobian $E/k$.  The number
of genus one curves $C'/k$ with $J(C) \cong J(C')$ which are $k$-isogenous
to $C$ is $N_C := \#(\Z/n\Z)^{\times}/(\pm 1)$, where $C \in H^1(K,E)$
has exact order $n$. In particular $N_C = 1$ if and only if
$n$ one of: $1, \ 2, \ 3,\ 4, \ 6$; and $N_C \ra \infty$ with $n$.
\end{cor}
\noindent
The proof is immediate from the previous proposition.  This result should be
compared with Amitsur's theorem: it is not true that two genus
one curves, even with common Jacobian, which have the same splitting fields
must be birational.
\begin{prop}
The field-isogeny class of a one-dimensional function field with respect to a number
field is finite.
\end{prop}
\noindent
Proof: When the genus is different from one, we have seen that
field-isogeny implies isomorphism, so it remains to look at the case
of $K$ a genus one function field with respect to a number field
$k$.  Fix some $k$-structure on $K$ (there are, of course, only
finitely ways to do this).  For the sake of clarity, let us
first show that there are only finitely many function fields
$K'/k$ which are isogenous to $K$ as $k$-algebras.  It will then
be easy to see that the proof actually gives finiteness of the
field-isogeny class.
\\ \indent
Let $C/k$ be the genus one curve such that $K = k(C)$; let
$K'/k$ be a function field such that there exists a homomorphism
$\iota: K' \ra K$, which on the geometric side corresponds to 
a finite morphism $\varphi: C \ra C'$, where $C'/k$ is another
genus one curve with $K' = k(C')$.  \\
\indent
By passing
to the Jacobian\footnote{Motivated by a desire to 
reverse as few arrows as possible, we choose to take the covariant Jacobian
functor, i.e., the Albanese rather than the Picard.}
of $\varphi/k$ we get an isogeny of elliptic curves
$J(C) \ra J(C')$.  By Shafarevich's theorem, the isogeny class of
an elliptic curve over a number field is finite , 
so it is enough to
bound the number of function fields $k(C')$ with a given Jacobian, say
$E'$.  Let $n$ be the common period of $C$ and $C'$, so that
$C' \in H^1(k,E')[n]$.  The set $S$ of places of $k$
containing the infinite places and all finite places $v$ such
that $C(k_v) = \emptyset$; this is a finite set.  But the existence
of $\varphi$ means that $C'(k_v) \neq \emptyset$ for all $v$ outside of
$S$, so that 
\[C' \in \ker(H^1(k,E)[n] \ra \prod_{v \not{\in} S} H^1(k_v,E)[n]). \]
But the finiteness of this kernel is extremely well-known, using
e.g. Hermite's discriminant bounds.  (Indeed, this is the key step
in the proof of the weak Mordell-Weil theorem; see e.g. [Silverman].)  
\\ \indent
Notice that we actually showed the following: a given genus
one function field $K/k$ dominates only finitely many other genus one
function fields.  (In fact $K$ dominates only finitely many function
fields in all, the genus zero case being taken care of by the
finiteness of the Brauer kernel $\kappa(C)$.)  
It follows that there are only finitely many isomorphism classes
of fields $L$ which admit a field isomorphic to $K$ as a finite extension;
this completes the proof.
\subsection{The proof of Theorem 10}
Let $k$ be a number field and $C_1/k$ a genus one curve of period $1, \ 2, \ 3, \ 4, \ $ or $6$ 
whose Jacobian $J(C_1)$ 
has no complex multiplication over $\overline{k}$ and is isolated in its 
isogeny class.  Let $K_1 = k(C_1)$ and $K_2$ be any finitely generated
field such that $K_1 \equiv K_2$.  By Pop's Theorem 3, $K_1$ and
$K_2$ are isogenous as fields, so $K_2$ is isomorphic as a field
to $k(C_2)$ where $C_2/k$ is another genus one curve.  By modifying
if necessary the $k$-structure on $C_2$, we get a finite morphism $\varphi: 
C_1 \ra C_2$ of $k$-schemes; passing to $J(\varphi)$ we deduce
that $J(C_1) \sim J(C_2)$ and hence by hypothesis that 
$E = J(C_1) \cong J(C_2)$.
By Proposition 19, $[C_2] = a[C_1]$ for some integer $a$.  Applying the same
argument with the roles of $C_1$ and $C_2$ interchanged, we get
that $[C_1]$ and $[C_2]$ generate the same cyclic subgroup of
$H^1(k,E)$, and by the hypothesis on the period of $C_1$ we conclude
$C_1 \cong C_2$.
\\ \\
We remark that with hypotheses as above but $C$ of arbitrary period 
$n$, we find that $k(C)$ could be elementarily equivalent only
to one of $\#(\Z/n\Z)^{\times}/(\pm 1)$ nonisomorphic function fields,
but distinguishing between these isogenous genus one curves with common
Jacobian seems quite difficult.
\\ \\
Finally, we must show that the assumption that $J(C)$ is isolated
can be removed at the cost of assuming the finiteness of the 
Mordell-Weil group $J(C)(k)$.  This is handled by the following
result, which is a modification of the (clever, and somewhat tricky) argument of 
[Pierce] to our arithmetic situation.
\begin{prop}
Let $K_1 = k(C_1)$ be the function field of a genus one curve over a number 
field.  Assume that $J(C_1)$ does not have complex multiplication (even) over
the algebraic closure of $k$, and that $J(C_1)(k)$ is finite.  Let
$K_2 \equiv K_1$ be any elementarily equivalent function field.  Then
$K_2 = k(C_2)$ is the function field of a genus one curve $C_2$ such
that $J(C_1) \cong J(C_2)$.
\end{prop}
\noindent
Proof: Let $K_2$ be a finitely generated function field such that 
$K_1 \equiv K_2$.  By Pop's Theorem A, we know
that $K_2$ is field-isogenous to $K_1$.  As above, this implies
the existence of $k$-structures on $K_1$ and $K_2$ such that
$K_1 = k(C_1), \ K_2 = k(C_2)$ and $\iota/k: C_1 \ra C_2$
is a finite morphism.  (Again we get by without using the
full strength of the notion of field-isogeny.)
\\ \\
Step 1: In search of a contradiction, we assume that the greatest common divisor 
of the degrees of all
finite maps $C_1 \ra C_2$ is divisible by some prime number $p$.  
\\ \\
Step 2: Because of Step 1 and the finitness of $J(C_2)(k)$, there is a 
finite list of \'etale maps $\lambda_i:
C_1 \ra C_i$ such that every map $C_1 \ra C_2$ factors through
some $\lambda_i$:
\[C_1 \stackrel{\lambda_i}{\ra} C_i \stackrel{\Psi_i}{\ra} C_2. \]
Step 3: We choose a smooth affine model for $C_2/k$ and let
$\overline{x} = (x_1,\ldots,x_n)$ denote coordinates.  The statement
``$\overline{x} \in C_2$'' can be viewed as first-order: let
$(P_j)$ be a finite set of generators for the ideal of $C_2$
in $k[\overline{x}]$; then $\overline{x} \in C_2$ is an abbreviation
for ``$\fa j \ P_j(\overline{x}) = 0$''.  For
each $i$, choose $b_i \in k(C_i)$ such that 
$k(C_i) = \Psi_i^*(k(C_2))(b_i)$, and let $g_i(X,Y) \in k[\overline{X},Y]$ be the
minimal polynomial for $b_i$ over $\Psi_i^*(k(C_2))$.  Finally, we
define a predicate $\overline{x} \in C_2 \wedge \neg \  \Con(\overline{x})$
with the meaning that $\overline{x}$ lies on $C_2$ and each coordinate
is not in $k$.  We must stress that this is to be regarded as a single
symbol -- we do not know how to define the constants in a function field
over a number field, but since $C_2$ by assumption has only finitely
many $k$-rational points, we can name them explicitly.  Consider the sentence:
\[ \fa \overline{x} \te y  \left(\overline{x} \in C_2 \wedge \neg \ \Con
(\overline{x})) \implies \bigvee_i 
g_i(\overline{x},y) = 0 \right) \]
Note well that $k(C_2)$ does not satisfy this sentence: take $\overline{x}$
to be any generic point of $C_2$.  But $k(C_1)$ does: giving such an element
$\overline{x} \in k(C_1)$ is equivalent to giving a field embedding
$\iota: k(C_2) \ra k(C_1)$, i.e. to a finite map $\iota: C_1 \ra C_2$.  
So $\iota = \Psi_i \circ
\lambda_i$ for some $i$, and we can take $y = \lambda_i^* b_i$:
\[g(\overline{x},y) = g(\iota^*\overline{a},\lambda^* b) = 
\lambda_i^* g(\Psi_i^* \overline{a},b) = 0, \]
with $\overline{x} = \iota^*(\overline{a})$, $\overline{a}$ a generic
point of $C_2$. So our sentence exhibits the elementary inequivalence of $k(C_1)$ and
$k(C_2)$, a contradiction.  
\\ \\
Step 4: Therefore the assumption of Step 1 is false, and it follows that there exist two
isogenies between the non-CM elliptic curves $C_1/\overline{k}$, 
$C_2/\overline{k}$ of coprime degree, and this easily implies that 
$j(C_1) = j(C_2)$.  In particular, the Jacobians $J(C_1)$ and $J(C_2)$
are isogenous elliptic curves with the same $j$-invariant and without
complex multiplication.  This implies that $J(C_1)$ and $J(C_2)$ are
isomorphic over $k$: indeed, let $\iota: J(C_1) \ra J(C_2)$ be any isogeny.
Then, $\iota/\overline{\Q}$ must have Galois group $\Z/n\Z \oplus \Z/n\Z$
for some $n$, so that $J(C_2) = J(C_1)/ker(\iota) = J(C_1)/J(C_1)[n] \cong J(C_1)$.
\\ \\
The end of the proof is the same as in the first case of the theorem: since
the period is $1, \ 2, \ 3, \ 4,$ or $6$, we may conclude $C_1 \cong C_2$.
\\ \\
Remark: The proof uses the finiteness of $J(C)(k)$ in two places: in order
to get around the nondefinability of $k$ in $K$, but also to get
around the fact that whereas over an algebraically closed field,
an arbitrary finite morphism of elliptic curves $E_1 \ra E_2$ can
be factored as $\iota \circ \tau$, over an arbitrary field
we can only claim the factorization $\tau \circ \iota$.

\end{document}